\documentclass[twoside,12pt]{article}

\usepackage{amsmath,amssymb,amsthm,mathrsfs}
\usepackage[top=2cm,bottom=2cm,left=2cm,right=2cm]{geometry}
\usepackage[hidelinks]{hyperref}

\newtheorem{theorem}{Theorem}[section]
\newtheorem{lemma}{Lemma}[section]
\newtheorem{conjecture}{Conjecture}

\newtheorem{remark}{Remark}[section]

\begin{document}

\title{\textbf{Asymptotic Bounds for $t(3,n)$ and an Application to $t(4,n)$}\footnote{The first author is supported by the Tianjin Municipal Education Commission Scientific Research Program Project (Grant No. 2025KJ133), the program of China Scholarship Council (Grant No. 202308120095) and Institute for Basic Science (IBS-R029-C4); the second author is supported by the National Natural Science Foundation of China (Nos. 12061059, 11601254, and 11551001) and the Qinghai Key Laboratory of Internet of Things Project (2017-ZJ-Y21).}}

\author{Meng Ji\footnote{School of Mathematical Sciences, and Institute of Mathematics and Interdisciplinary Sciences,Tianjin Normal University, Tianjin, China. {\tt
mji@tjnu.edu.cn}} \footnote{Extremal Combinatorics and Probability Group (ECOPRO)
 Institute for Basic Science (IBS), Daejeon, South Korea}, \ \ Yaping Mao\footnote{School of Mathematics
and Statistis, Qinghai Normal University, Xining, Qinghai 810008,
China. {\tt maoyaping@ymail.com}}, \ \ Ingo
Schiermeyer\footnote{Technische Universit{\"a}t Bergakademie
Freiberg, Institut f{\"u}r Diskrete Mathematik und Algebra, 09596
Freiberg, Germany. {\tt Ingo.Schiermeyer@math.tu-freiberg.de}}}
\date{}
\maketitle

\begin{abstract}
Let $X\subseteq V(G)$ be a vertex subset of a simple graph $G$. The set $X$ is called \emph{irredundant} if every vertex $x\in X$ is either isolated in the induced subgraph $G[X]$ or possesses a private neighbor $y\in V(G)\setminus X$ that is adjacent to $x$ and to no other vertex of $X$.  
The \emph{mixed Ramsey number} $t(m,n)$ is defined as the smallest integer $N$ such that every red-blue coloring of the edges of the complete graph $K_N$ contains either an $m$-element irredundant set in the blue graph or an $n$-element independent set in the red graph. The \emph{irredundant Ramsey number} $s(m,n)$ is defined analogously, with the independent set replaced by an $n$-element irredundant set in the red graph.  

In this paper, we prove that for every $\varepsilon>0$, the upper bound
\[
        t(3,n)\le (2+\varepsilon)\frac{n^{3/2}}{\sqrt{\log n}}
\]
holds for all sufficiently large $n$. When combined with a lower bound due to Krivelevich, this result yields asymptotic behavior of both $s(3,n)$ and $t(3,n)$. As an application, we employ the recent lower bound $r(4,n)=\Omega\bigl(n^3/(\log n)^4\bigr)$ established by Mattheus and Verstraete to show that
\[
        \lim_{n\to\infty}\frac{t(4,n)}{r(4,n)}=0.
\]

\vspace{2mm}

\end{abstract}

\section{Introduction}

The Ramsey number $r(G,H)$ is the minimum integer $N$ such that every
red-blue coloring of the edges of the complete graph $K_N$ contains a blue
copy of $G$ or a red copy of $H$.  In this paper, we consider a classical
variant of Ramsey theory involving irredundant sets.

In 1978, Cockayne, Hedetniemi, and Miller \cite{Cockayne-Hedetniemi-Miller}
introduced the concept of irredundance in connection with domination in graphs.
Let $G=(V,E)$ be a simple graph.  A set of vertices $X\subseteq V$ is
\emph{irredundant} if each vertex $x\in X$ is either isolated in the induced
subgraph $G[X]$ or has a \emph{private neighbor} $y\in V\setminus X$ that is
adjacent to $x$ and to no other vertex of $X$.

The \emph{irredundant Ramsey number} $s(m,n)$ is the minimum integer $N$ such
that every red-blue coloring of the edges of $K_N$ contains an $m$-element
irredundant set in the blue graph or an $n$-element irredundant set in the red
graph.  The \emph{mixed Ramsey number} $t(m,n)$ is the minimum integer $N$ such
that every red-blue coloring of the edges of $K_N$ contains an $m$-element
irredundant set in the blue graph or an $n$-element independent set in the red
graph.  Since every independent set is irredundant, we have
\[
        s(m,n)\le t(m,n)\le r(m,n).
\]

Brewster, Cockayne, and Mynhardt \cite{BrewsterCockayneMynhardt} proved that,
in a red-blue coloring of the edges of a complete graph, the blue graph contains
a $3$-element irredundant set if and only if the red graph contains a triangle
or an induced cycle of length six.  This characterization is the starting point
of our proof.

Chen, Hattingh, and Rousseau \cite{chen1993}, Erd\H{o}s and Hattingh
\cite{ErdosHattingh}, and Krivelevich \cite{Krivelevich} obtained several
asymptotic bounds for irredundant Ramsey numbers and mixed Ramsey numbers.
In particular, Krivelevich proved the lower bound
\[
        s(3,n)\ge c\left(\frac{n}{\log n}\right)^{5/4}
\]
for some positive constant $c$.  Rousseau and Speed \cite{RousseauSpeed}
showed that
\[
        t(3,n)\le \frac{5n^{3/2}}{\sqrt{\log n}}
\]
for all sufficiently large $n$.  We give a corrected Shearer-Hattingh type
proof yielding the following slightly sharper asymptotic constant.

\begin{theorem}\label{th2-3}
There is a positive constant $c_1$ such that, for every $\varepsilon>0$, there
exists $n_0=n_0(\varepsilon)$ with
\[
        c_1\left(\frac{n}{\log n}\right)^{5/4}
        < s(3,n)\le t(3,n)
        \le (2+\varepsilon)\frac{n^{3/2}}{\sqrt{\log n}}
\]
for all $n\ge n_0$.
\end{theorem}

Chen, Hattingh, and Rousseau \cite{chen1993} also proposed the following
conjecture; see also Mynhardt and Roux \cite{Mynhardt-Roux}.

\begin{conjecture}[Chen, Hattingh, and Rousseau]\label{conjecture}
For each fixed $m\ge 4$,
\[
        \lim_{n\to\infty}\frac{t(m,n)}{r(m,n)}=0.
\]
\end{conjecture}

Using Theorem~\ref{th2-3} and the recent lower bound on $r(4,n)$ due to
Mattheus and Verstraete \cite{MattheusVerstraete}, we verify the conjecture for
$m=4$.

\begin{theorem}\label{them12}
\[
        \lim_{n\to\infty}\frac{t(4,n)}{r(4,n)}=0.
\]
\end{theorem}

All logarithms are natural. For a graph $G$, we write $V(G)$ and
$E(G)$ for its vertex set and edge set.  For a vertex $v\in V(G)$, let $N(v)$
denote the open neighborhood of $v$ and $N[v]=N(v)\cup\{v\}$ its closed
neighborhood.  If distances are taken in a graph $G$, we write
\[
        D_i(v)=\{u\in V(G): d_G(u,v)=i\},\qquad
        D_{>i}(v)=\{u\in V(G): d_G(u,v)>i\}.
\]
Vertices not connected to $v$ in $G$ are regarded as having distance larger
than every integer.  In the proof below, these distance sets are used in the
red graph unless otherwise specified.

\section{The upper bound for \texorpdfstring{$t(3,n)$}{t(3,n)}}

We use three standard ingredients.  The first is a structure theorem of
Hattingh.

\begin{theorem}[Hattingh \cite{Hattingh}]\label{Hattingh}
Let $(R,B)$ be a red-blue coloring of the edges of a complete graph $K_N$ in
which the blue graph $B$ contains no $3$-element irredundant set.  Fix a vertex
$v$.  Partition $N_B(v)$ into two parts $D_2(v)$ and $D_{>2}(v)$, where the
distance is measured in the red graph $R$.  Thus a vertex $u\in N_B(v)$ belongs
to $D_2(v)$ if $d_R(u,v)=2$, and belongs to $D_{>2}(v)$ otherwise
(that is, its red distance from $v$ is greater than two, possibly infinite).
If $X\subseteq N_B(v)$ contains at most one vertex from $D_{>2}(v)$, then the red
induced graph $R[X]$ is bipartite.
\end{theorem}

The second ingredient is the following characterization of blue irredundant
sets of size three.

\begin{theorem}[Brewster, Cockayne, and Mynhardt \cite{BrewsterCockayneMynhardt}]\label{lem2}
In a red-blue coloring of the edges of a complete graph, the blue graph contains
a $3$-element irredundant set if and only if the red graph contains a triangle
or an induced cycle of length six.
\end{theorem}

The third ingredient is Shearer's lower bound for the independence number of a
triangle-free graph.

\begin{lemma}[Shearer \cite{shearer1983}]\label{shearer}
Let
\[
        f(x)=\frac{x\log x-(x-1)}{(x-1)^2},\qquad f(1)=\frac12.
\]
If $G$ is a triangle-free graph of order $N$ and average degree $\bar d$, then
\[
        \alpha(G)\ge Nf(\bar d).
\]
Moreover, $f$ is decreasing on $[1,\infty)$.
\end{lemma}
For completeness, the monotonicity follows from
\[
        f'(x)=\frac{2x-2-(x+1)\log x}{(x-1)^3}\le 0\qquad (x>1),
\]
since $(x+1)\log x-2x+2$ is nonnegative on $[1,\infty)$.

\begin{proof}[Proof of Theorem~\ref{th2-3}]
The lower bound is due to Krivelevich \cite{Krivelevich}, and the inequality
$s(3,n)\le t(3,n)$ follows because every independent set is irredundant.  It
remains to prove the asserted upper bound for $t(3,n)$.

Fix $\varepsilon>0$.  Put
\[
        \eta=\frac{\varepsilon}{2},\qquad C=2+\eta.
\]
Define $F(1)=0$, and for $x>1$ define
\[
        F(x)=C\frac{x^{3/2}}{\sqrt{\log x}}.
\]
We first prove an auxiliary bound with an additive constant:
\[
        t(3,n)\le \lceil F(n)+L\rceil
\]
for all sufficiently large $n$, where $L$ is a fixed constant depending only on
$\varepsilon$.

Let
\[
        q=q(n)=\sqrt{n\log n}.
\]
Since
\[
        f(q)=\frac{q\log q-(q-1)}{(q-1)^2}
             =(1+o(1))\frac{\log q}{q}
             =\left(\frac12+o(1)\right)\sqrt{\frac{\log n}{n}},
\]
we have
\[
        F(n)f(q)=\left(\frac C2+o(1)\right)n.
\]
As $C>2$, it follows that
\begin{equation}\label{case1-prep}
        F(n)f\bigl(\sqrt{n\log n}\bigr)>n
\end{equation}
for all sufficiently large $n$.

Also, for $x>1$,
\[
        F'(x)=C\frac{x^{1/2}}{\sqrt{\log x}}
        \left(\frac32-\frac{1}{2\log x}\right).
\]
Uniformly for $x\in[n-q,n]$, we have
\[
        F'(x)=\left(\frac{3C}{2}+o(1)\right)
        \frac{n^{1/2}}{\sqrt{\log n}}.
\]
Therefore
\[
        F(n)-F(n-q)=\int_{n-q}^{n}F'(x)\,dx
        =\left(\frac{3C}{2}+o(1)\right)n.
\]
Since $C>2$, we get
\begin{equation}\label{case2-prep}
        F(n)-F\bigl(n-\sqrt{n\log n}\bigr)>3n+1
\end{equation}
for all sufficiently large $n$.

Choose $n_1$ sufficiently large so that \eqref{case1-prep} and
\eqref{case2-prep} hold for every $n\ge n_1$, so that $F(n)>2n$ for every
$n\ge n_1$, and so that $n-\sqrt{n\log n}\ge2$ for every $n\ge n_1$.
Note also that $F$ is increasing on $[2,\infty)$.  Choose $L>0$ so large that
\begin{equation}\label{small-k}
        t(3,k)\le \lceil F(k)+L\rceil
\end{equation}
for every integer $1\le k<n_1$.  This is possible because there are only
finitely many such $k$.

We now prove by induction on $n\ge n_1$ that
\begin{equation}\label{induction-claim}
        t(3,n)\le \lceil F(n)+L\rceil.
\end{equation}
Let
\[
        N=\lceil F(n)+L\rceil,
\]
and consider an arbitrary red-blue coloring of the edges of $K_N$.  Let $R$ be
the red graph and $B$ the blue graph.  Suppose, for a contradiction, that $B$
contains no $3$-element irredundant set and that $R$ contains no independent set
of size $n$.  By Theorem~\ref{lem2}, the red graph $R$ is triangle-free and has
no induced cycle of length six.  Let $\bar d_R$ be the average degree of $R$.

\medskip
\noindent\textbf{Case 1.} $\bar d_R<q$.

If $\bar d_R<1$, then the elementary bound $\alpha(R)\ge N/(\bar d_R+1)$ gives
$\alpha(R)>N/2\ge F(n)/2>n$, by the choice of $n_1$.  This is a contradiction.
Thus we may assume $\bar d_R\ge1$.  Since $f$ is decreasing on $[1,\infty)$,
Shearer's lemma gives
\[
        \alpha(R)\ge Nf(\bar d_R)\ge Nf(q)
        \ge (F(n)+L)f(q)
        \ge F(n)f(q)>n,
\]
where the last inequality is \eqref{case1-prep}.  This contradicts the
assumption that $R$ has no independent set of size $n$.

\medskip
\noindent\textbf{Case 2.} $\bar d_R\ge q$.

Choose a vertex $v$ with
\[
        d_R(v)\ge\bar d_R\ge q.
\]
Since $R$ is triangle-free, $N_R(v)$ is an independent set in $R$.  Hence
\[
        d_R(v)=|N_R(v)|\le n-1.
\]
Put
\[
        d=d_R(v),\qquad k=n-d.
\]
Then $1\le k<n$.  Here $D_2(v)$ and $D_{>2}(v)$ denote the two parts of
$N_B(v)$ from Theorem~\ref{Hattingh}, with distances measured in the red graph
$R$.

Theorem~\ref{Hattingh}, applied with $X=D_2(v)$, implies that $R[D_2(v)]$ is
bipartite.  If one part of this bipartite graph had at least $n$ vertices, then
$R$ would contain an independent set of size $n$, contrary to our assumption.
Thus
\begin{equation}\label{D2-bound}
        |D_2(v)|\le 2n-2.
\end{equation}

There are no red edges between $N_R(v)$ and $D_{>2}(v)$; otherwise a vertex of
$D_{>2}(v)$ would be at red distance two from $v$.  Therefore every independent
set in $R[D_{>2}(v)]$, together with $N_R(v)$, is an independent set in $R$.
Consequently
\[
        \alpha(R[D_{>2}(v)])+d\le \alpha(R)<n,
\]
and so
\[
        \alpha(R[D_{>2}(v)])<n-d=k.
\]
The induced blue graph on $D_{>2}(v)$ also contains no $3$-element irredundant
set.  Hence, by the definition of $t(3,k)$,
\begin{equation}\label{Dgreater-bound}
        |D_{>2}(v)|<t(3,k).
\end{equation}

Since
\[
        V(K_N)=\{v\}\cup N_R(v)\cup D_2(v)\cup D_{>2}(v),
\]
we obtain from \eqref{D2-bound} and \eqref{Dgreater-bound} that
\[
        N=1+d+|D_2(v)|+|D_{>2}(v)|
        <1+d+(2n-2)+t(3,k)<t(3,k)+3n.
\]
If $k<n_1$, then \eqref{small-k} gives
$t(3,k)\le\lceil F(k)+L\rceil$; if $k\ge n_1$, the same bound follows from the
induction hypothesis.  Since $d\ge q$, we have $k=n-d\le n-q$.  If $k=1$, then
$F(k)=0\le F(n-q)$; if $k\ge2$, the same conclusion follows from the
monotonicity of $F$ on $[2,\infty)$.  Thus
\[
        N<\lceil F(k)+L\rceil+3n
        \le F(k)+L+1+3n
        \le F(n-q)+L+1+3n.
\]
By \eqref{case2-prep},
\[
        F(n-q)+L+1+3n<F(n)+L.
\]
Therefore
\[
        N<F(n)+L\le\lceil F(n)+L\rceil=N,
\]
a contradiction.

Both cases give contradictions.  Hence every red-blue coloring of the edges of
$K_N$ contains either a $3$-element irredundant set in the blue graph or an
$n$-element independent set in the red graph.  This proves
\eqref{induction-claim} for every $n\ge n_1$.

Finally, since $L+1=o(n^{3/2}/\sqrt{\log n})$, for all sufficiently large $n$ we
have
\[
        \lceil F(n)+L\rceil
        \le (2+\varepsilon)\frac{n^{3/2}}{\sqrt{\log n}}.
\]
This completes the proof of Theorem~\ref{th2-3}.
\end{proof}

\section{Application to \texorpdfstring{$t(4,n)$}{t(4,n)}}

We first recall the standard recurrence for mixed Ramsey numbers.

\begin{lemma}\label{recurrence}
For all $m,n\ge3$,
\[
        t(m,n)\le t(m-1,n)+t(m,n-1).
\]
\end{lemma}

\begin{proof}
Let $N=t(m-1,n)+t(m,n-1)$, and consider a red-blue coloring of the edges of
$K_N$.  Fix a vertex $v$.  If
\[
        d_R(v)\ge t(m-1,n),
\]
then the coloring induced on $N_R(v)$ contains either an $(m-1)$-element
irredundant set in the blue graph or an $n$-element independent set in the red
graph.  In the first case, adding $v$ gives an $m$-element irredundant set in
the blue graph.  Indeed, $v$ is blue-isolated from $N_R(v)$, and any private
neighbor that certifies irredundance inside $N_R(v)$ remains private after $v$
is added because all vertices of $N_R(v)$ are red-adjacent to $v$.  In the
second case the conclusion follows immediately.

If $d_R(v)<t(m-1,n)$, then
\[
        d_B(v)=N-1-d_R(v)\ge t(m,n-1).
\]
The coloring induced on $N_B(v)$ contains either an $m$-element irredundant set
in the blue graph or an $(n-1)$-element independent set in the red graph.  In
the first case we are done; in the second case, adding $v$ gives an
$n$-element independent set in the red graph, because $v$ has no red edges to
$N_B(v)$.  This proves the recurrence.
\end{proof}

The following lower bound for $r(4,n)$ is due to Mattheus and Verstraete.

\begin{theorem}[Mattheus and Verstraete \cite{MattheusVerstraete}]\label{MV}
There exists a positive constant $c$ such that
\[
        r(4,n)\ge c\frac{n^3}{(\log n)^4}
\]
for all sufficiently large $n$.
\end{theorem}

\begin{proof}[Proof of Theorem~\ref{them12}]
By Lemma~\ref{recurrence},
\[
        t(4,n)\le t(4,n-1)+t(3,n).
\]
Iterating this inequality from $3$ to $n$ and absorbing the fixed term
$t(4,2)$, as well as the finitely many initial $t(3,k)$ terms before
Theorem~\ref{th2-3} applies, into the implicit constant, we obtain
\[
        t(4,n)
        =O\left(\sum_{k=3}^{n}\frac{k^{3/2}}{\sqrt{\log k}}\right).
\]
For the sum, the terms with $k\le\sqrt n$ contribute at most
$O(n^{5/4})$, while for $k>\sqrt n$ we have
$\sqrt{\log k}\ge\sqrt{(\log n)/2}$.  Therefore
\[
        \sum_{k=3}^{n}\frac{k^{3/2}}{\sqrt{\log k}}
        =O\left(n^{5/4}\right)
         +O\left(\frac{1}{\sqrt{\log n}}\sum_{k\le n}k^{3/2}\right)
        =O\left(\frac{n^{5/2}}{\sqrt{\log n}}\right).
\]
Hence
\[
        t(4,n)=O\left(\frac{n^{5/2}}{\sqrt{\log n}}\right).
\]
By Theorem~\ref{MV},
\[
        r(4,n)\ge c\frac{n^3}{(\log n)^4}
\]
for all sufficiently large $n$.  Therefore
\[
        0\le\frac{t(4,n)}{r(4,n)}
        \le O\left(
        \frac{n^{5/2}/\sqrt{\log n}}{n^3/(\log n)^4}
        \right)
        =O\left(\frac{(\log n)^{7/2}}{\sqrt n}\right)\to0.
\]
This proves the theorem.
\end{proof}

\begin{remark}
The proof above does not establish an upper bound of order
$n^{5/4}/\log n$ for $t(3,n)$.  In the Shearer-Hattingh framework used here,
the recurrence obtained in the high average degree case has an additive term of
order $n$, and this naturally balances with the low average degree case at the
scale $n^{3/2}/\sqrt{\log n}$.
\end{remark}

\vspace{4mm}
\noindent{\bf Acknowledgment.} We thank Prof. Guantao Chen for reading some
parts of the paper and for helpful comments, and Prof. Xian'an Jin for his valuable suggestions on the exposition of this article.

\end{document}